\documentclass{amsart}
\usepackage{graphicx}
\usepackage{amssymb}
\usepackage{amsfonts}
\swapnumbers
\newtheorem{theorem}{Theorem}[section]

\theoremstyle{definition}

\newtheorem{remark}[theorem]{Remark}

\numberwithin{equation}{section}
 \theoremstyle{plain}    
 
 \numberwithin{equation}{section} %% Comment out for sequentially-numbered
 \numberwithin{figure}{section} %% Comment out for sequentially-numbered
 \theoremstyle{plain}    
 \theoremstyle{plain}    
 \theoremstyle{remark}    
 \newtheorem*{acknowledgement*}{Acknowledgement} 
\newcommand{\cB}{{\mathcal B}}

\newcommand{\cF}{{\mathcal F}}

\newcommand{\cJ}{{\mathcal J}}
\newcommand{\cK}{{\mathcal K}}

\newcommand{\Om}{{\Omega}}
\newcommand{\om}{{\omega}}
\newcommand{\ve}{{\varepsilon}}
\newcommand{\del}{{\delta}}

\newcommand{\gam}{{\gamma}}
\newcommand{\Gam}{{\Gamma}}

\newcommand{\sig}{{\sigma}}
\newcommand{\al}{{\alpha}}
\newcommand{\be}{{\beta}}

\newcommand{\la}{{\lambda}}

\newcommand{\bbZ}{{\mathbb Z}}
\newcommand{\bbI}{{\mathbb I}}
\begin{document}
\title[]{Nonconventional Poisson limit theorems}%
 \vskip 0.1cm 
 \author{ Yuri Kifer\\
\vskip 0.1cm
 Institute  of Mathematics\\
Hebrew University\\
Jerusalem, Israel}%
\address{
Institute of Mathematics, The Hebrew University, Jerusalem 91904, Israel}
\email{ kifer@math.huji.ac.il}%

\thanks{ }
\dedicatory{  }
 \date{\today}
\begin{abstract}\noindent
The classical Poisson theorem says that if $\xi_1,\xi_2,...$ are i.i.d.
0--1 Bernoulli random variables taking on 1 with probability $p_n\equiv
\la/n$ then the sum $S_n=\sum_{i=1}^n\xi_i$ is asymptotically in $n$ 
Poisson distributed with the parameter $\la$. It turns out that this result
can be extended to sums of the form $S_n=\sum_{i=1}^n\xi_{q_1(i)}\cdots
\xi_{q_\ell(i)}$ where now $p_n\equiv(\la/n)^{1/\ell}$ and $1\leq q_1(i)
<\cdots <q_\ell(i)$ are integer valued increasing functions. We obtain also 
Poissonian limit for numbers of arrivals to small sets of $\ell$-tuples 
$X_{q_1(i)},...,X_{q_\ell(i)}$ for some Markov chains $X_n$ and for numbers
of arrivals of $T^{q_1(i)}x,...,T^{q_\ell(i)}x$ to small cylinder sets for 
typical points $x$ of a subshift of finite type $T$. 
\end{abstract}
%\footnotetext[1]{} 
\maketitle
\markboth{Yu.Kifer}{Nonconventional averaging} 
\renewcommand{\theequation}{\arabic{section}.\arabic{equation}}
\pagenumbering{arabic}

\section{Introduction}\label{sec1}\setcounter{equation}{0}

The classical Poisson limit theorem taught in the first probability 
course says that if
\begin{equation}\label{1.1}
\lim_{n\to\infty}np_n=\la>0,\,\, p_n>0
\end{equation}
then the binomial distribution with parameters $(n,p_n)$ converges as
$n\to\infty$ to the Poisson distribution with parameter $\la$. In other
 words, if $\xi_1^{(n)},\xi_2^{(n)},\xi_3^{(n)}...,\, n=1,2,...$ is an
 array of independent Bernoulli random variables satisfying
 \begin{equation}\label{1.2}
 p_n=P\{\xi_i^{(n)}=1\}=1-P\{\xi_i^{(n)}=0\}
 \end{equation}
 and (\ref{1.1}) holds true then the sum
 \begin{equation}\label{1.3}
 S_n=\sum_{l=1}^n\xi_l^{(n)}
 \end{equation}
 converges in distribution to a Poisson random variable with the parameter 
 $\la$.
 
 It turns out that assuming  
 \begin{equation}\label{1.4}
 \lim_{n\to\infty}np_n^\ell=\la
 \end{equation}
 the above result can be extended to "nonconventional" sums of the
 form
 \begin{equation}\label{1.5}
 S_n=\sum_{l=1}^n\xi_{q_1(l)}^{(n)}\xi_{q_2(l)}^{(n)}\cdots
 \xi_{q_\ell(l)}^{(n)}
 \end{equation}
 where $l\leq q_1(l)<q_2(l)<...<q_\ell(l)$ are increasing functions taking on 
 integer values on integers. The name "nonconventional" comes from \cite{Fu}
where ergodic theorems for sums of the form (\ref{1.5}) were studied. Recently,
 strong laws of large numbers and central limit theorems type results were
 obtained in \cite{Ki} and \cite{KV} even for more general expressions.
 
 We will consider also a nonconventional Poisson limit theorem for sums of
 the form
 \begin{equation}\label{1.6}
 S_n=\sum_{l=1}^n\prod_{j=1}^\ell\bbI_{\Gam_n}(X_{q_j(l)})
 \end{equation}
 (where $\bbI_\Gam$ is the indicator of a set $\Gam$), which counts the 
 number of arrivals by a Markov chain $X_0,X_1,...$ to small sets $\Gam_n$
  at all times $q_j(l),\, j=1,...,\ell$ when $l$ runs from 1 to
 $n$. We suppose that the Markov chain has bounded transition densities which
 satisfy a Doeblin type condition. This ensures existence of a unique
 invariant (probability) measure $\mu$. Assuming that
 \begin{equation}\label{1.7}
 \lim_{n\to\infty}n(\mu(\Gam_n))^\ell=\la>0
 \end{equation}
 we show that $S_n$ converges in distribution to a Poisson random variable
 provided $q_{i+1}(l)-q_i(l)\to\infty,\, l=1,...,\ell-1$ as $l\to\infty$.
 
 Another nonconventional Poisson limit theorem we deal with in this paper 
 concerns multiple arrivals to shrinking cylinder sets by subshifts. Namely,
 we show that the expressions of the form 
 \begin{equation}\label{1.8}
 S_n=\sum_{l=1}^{\la/(P(B_n))^\ell}\prod_{j=1}^\ell\bbI_{B_n}\circ 
 T^{q_j(l)}
 \end{equation}
 have almost surely asymptotically in $n$ the Poisson distribution with a
 parameter $\la$ provided
 \begin{equation}\label{1.9}
 \lim_{n\to\infty}n(P(B_n))^\ell=\la
 \end{equation}
 where $B_n$'s are certain cylinder sets of the subshift
 space $\Om$, $P$ is a Gibbs measure, $T$ is a left shift and $l\leq q_j(l),
 \,j=1,...,\ell$ are increasing integer valued functions with sufficiently 
 fast growth of differences $q_{j+1}(l)-q_j(l)$ as $l\to\infty$. This 
 assertion generalizes in the nonconventional direction the results from
 \cite{Pi} and \cite{De}. We mention also a related result in the 
 conventional setup about Poisson limits for numbers of arrivals to small
  shrinking sets by hyperbolic dynamical systems obtained in \cite{Do}.
 
 Observe that if
 \[
 \tau_{B_n}(\om)=\min\{ l:\, T^{q_j(l)}\om\in B_n,\,\forall j=1,...,\ell\}
 \]
 then the above result yields that asymptotically the distribution of 
 $(P(B_n))^\ell\tau_{B_n}$ is exponential since
 \[
 P\{(P(B_n))^\ell\tau_{B_n}>\la\}=P\{ S_n=0\}\to e^{-\la}\,\,\mbox{as}\,\,
 n\to\infty.
 \]
 More advanced results about limiting exponential distributions of
 properly normalized first return times to small shrinking sets 
 (see, for instance, \cite{Ab}, \cite{AS} and references there) and about
 Poissonian asymptotical behavior of distributions of numbers of arrivals
 to small shrinking sets (see, for instance, \cite{HV1}, \cite{HV2} and 
 references there) should be possible to derive in the nonconventional 
 framework, as well.

  \section{Preliminaries and main results}\label{sec2}\setcounter{equation}{0}

We start with a probability space $(\Om,\cF,P)$ and an array of independent
Bernoulli 0--1 random variables $\xi_1^{(n)},\xi_2^{(n)},...,\, n=1,2,...$
taking on value 1 with probability $p_n\in(0,1)$. Our setup includes also
increasing functions $l\leq q_1(l)<q_2(l)<...<q_\ell(l)$ taking on integer
values on the integers. For any set $\Gam\subset\bbZ_+$ of nonnegative integers
put
\begin{equation}\label{2.1}
P_\la(\Gam)=\sum_{l\in\Gam}e^{-\la}\frac {\la^l}{l!}
\end{equation}
which is the probability assigned to $\Gam$ by the Poisson distribution with
a parameter $\la>0$.

\begin{theorem}\label{thm2.1} Set $\la_n=np^\ell_n$. Then
\begin{equation}\label{2.2}
\sup_{\Gam\subset\bbZ_+}|P\{ S_n\in\Gam\}-P_{\la}(\Gam)|\leq (2\ell^2+1)p_n+
2|\la-\la_n|e^{\max(\la,\la_n)},
\end{equation}
where $S_n$ is defined by (\ref{1.5}), and the right hand side of (\ref{2.2})
tends to zero provided (\ref{1.4}) holds true.
\end{theorem}

We will prove this assertion in Section \ref{sec3} relying on Poisson 
approximation results for families of dissociated random variables from
\cite{BE} and \cite{BHJ}. We will show also that the convergence of
$P\{ S_n\in\Gam\}$ to $P_\la(\Gam)$ can be derived relying on the result
of \cite{Se}, as well, though this does not give the speed of convergence
as in (\ref{2.2}).

As a natural intermediate step from the independent case to a stationary 
(dynamical systems) case we exhibit a nonconventional Poisson limit theorem
for Markov chains. Namely, let $X_0,X_1,...$ be a Markov chain on a measurable
state space $(M,\cB)$ whose one-step and $n$-step transition probabilities 
$P(x,\cdot)$ and $P(n,x,\cdot)$ satisfy the conditions
\begin{equation}\label{2.3}
P(x,\Gam)\leq Cm(\Gam)\,\,\,\mbox{and}\,\,\, P(n_0,x,\Gam)\geq C^{-1}m(\Gam)
\,\,\,\forall\Gam\subset M,\,\Gam\in\cB
\end{equation}
for some probability measure $m$ on $(M,\cB)$, an integer $n_0\geq 1$ and a 
constant $C>0$.
Denote by $P_x,\, x\in M$ the probability on the path space provided $X_0=x$
and let $\nu$ be an arbitrary initial distribution (a probability measure
on $M$). The Markov chain $X_0,X_1,X_2,...$ will be considered now with
respect to the probability $P_\nu=\int_MP_xd\nu(x)$ on the path space.
Assume also that
$l\leq q_1(l)<q_2(l)<...<q_\ell(l)$ is a sequence of integer valued increasing
functions such that 
\begin{equation}\label{2.4}
\lim_{l\to\infty}(q_{i+1}(l)-q_i(l))=\infty\,\,\,\forall i=1,2,...,\ell -1.
\end{equation}
Let $\mu$ be the unique invariant measure of the Markov chain above which
exists in view of the (strong) Doeblin condition (\ref{2.3}) (see, for
instance, \cite{Doo}, \S 5, Ch.V).
\begin{theorem}\label{thm2.2}
Let $\Gam_n\in\cB$ be a sequence of measurable subsets of $M$ such that 
(\ref{1.7}) holds true and let $P_\nu$ be a probability on the path space 
$\Om$ of the Markov chain corresponding to any initial distribution $\nu$. 
Then $S_n$ 
defined by (\ref{1.6}) converges in distribution as $n\to\infty$ on the 
probability space $(\Om,P)$ to a Poisson random variable with the 
parameter $\la$.
\end{theorem}

Next, we consider another setup where $\Om$ is a space of sequences determined
by a 0--1 matrix $A=(\al_{ij},\, 1\leq i,j\leq\iota)$, namely
\begin{equation}\label{2.5}
\Om=\{\om=(\om_0,\om_1,\om_2,...):\, 1\leq\om_i\leq\iota\,\,\mbox{and}\,\,
\al_{\om_i\om_{i+1}}=1\,\,\forall i\geq 0\}.
\end{equation}
The space $\Om$ together with the left shift $T$ acting by $(T\om)_i=\om_{i+1}$
is called a subshift of finite type (\cite{Bo}). We assume that $A^\wp$ is
a positive matrix for some $\wp>0$ which makes $T$ topologically mixing. Let
$\phi$ be a H\" older continuous function on $\Om$ with respect to the
metric 
\[
d(\om,\tilde\om)=\exp(-\min\{ i\geq 0:\,\om_i\ne\tilde\om_i\}).
\]
There exists a unique $T$-invariant Gibbs probability $P$ corresponding to 
$\phi$ (see \cite{Bo}) characterized by the property that for some $\Pi$
(called the topological pressure), any cylinder set
$[a_0,a_1,...,a_{n-1}]=\{\om\in\Om:\, \om_i=a_i,\, i=0,1,...,n-1\}$
and each $\om\in[a_0,a_1,...,a_{n-1}]$,
\begin{equation}\label{2.6}
C^{-1}\leq\frac {P([a_0,a_1,...,a_{n-1}])}{\exp(-\Pi n+\sum_{i=0}^{n-1}
\phi\circ T^i(\om))}\leq C
\end{equation}
where $C>0$ depends only on $\phi$. Denote by $\cF_n$ the (finite)
$\sig$-algebra generated by all cylinder sets $[a_0,a_1,...,a_{n-1}]\subset
\Om$ and for each $\om^*=(\om^*_0,\om^*_1,...)\in\Om$ set $C_n(\om^*)=
[\om_0^*,\om_1^*,...,\om_{n-1}^*]$. Next, let 
$l\leq q_1(l)<q_2(l)<...<q_\ell(l)$ be a sequence
 of increasing functions taking on integer values on integers and such that
 for some $c,\gam>0$ and all $l\geq 1$,
 \begin{equation}\label{2.7}
 q_{i+1}(l)-q_i(l)\geq c(\ln l)^{1+\gam},\, i=1,2,...,\ell-1.
 \end{equation}

\begin{theorem}\label{thm2.3} For some $s\geq 0$ and each $\om^*\in\Om$ 
let $B_n(\om^*)\subset C_n(\om^*)$, $B_n(\om^*)\in\cF_{n+[s\ln n]},\, 
n=1,2,...$ be arbitrary sequences of sets. For each sequence $N_n(\om^*),
\, n=1,2,...$ satisfying
\begin{equation}\label{2.8}
\lim_{n\to\infty}N_n(\om^*)\big(P(B_n(\om^*))\big)^\ell=\la>0
\end{equation}
set
\[
S_{n,\om^*}(\om)=\sum_{l=0}^{N_n(\om^*)}\prod_{j=1}^\ell\bbI_{B_n(\om^*)}
\circ T^{q_j(l)}(\om).
\]
Then for $P$-almost all $\om^*\in\Om$,
\begin{equation}\label{2.9}
\lim_{n\to\infty}P\{\om\in\Om:\, S_{n,\om^*}(\om)=k\}=
e^{-\la}\frac {\la^k}{k!}.
\end{equation}
\end{theorem}

We will prove this assertion modifying appropriately the technique from
\cite{Pi} and \cite{De} while relying on the basic result from \cite{Se} which
became a major tool for deriving Poissonian type limit theorems in dynamical
systems.

\begin{remark}\label{rem2.4}
In Theorems \ref{thm2.2} and \ref{thm2.3} we rely on \cite{Se} which does not
provide speed of convergence to the Poisson distribution as in Theorem 
\ref{thm2.1}. Still, verifying conditions of \cite{Se} we obtain certain 
estimates for speed of convergence of relevant quantities there so relying on
the quantitive version of \cite{Se} obtained in \cite{HV1} we can obtain
some estimates on speed of convergence in Theorems \ref{thm2.2} and 
\ref{thm2.3}, as well.
\end{remark}

\begin{remark}\label{rem2.5}
It is well known that already in the conventional $\ell=1$ setup the assertion
of Theorem \ref{thm2.3} holds true not for all $\om^*$ but only for almost
all $\om^*$ and a corresponding example appears already in \cite{Pi}. In fact,
when $\om^*$ is a periodic point then already in the conventional $\ell=1$ 
setup the distribution of $S_{n,\om^*}$ from Theorem \ref{thm2.3} will converge 
to a compound Poisson distribution (see \cite{HV2}) and not to a Poisson one.
The nonconvergence to a Poisson distribution is usually easy to see in this
case checking that second moments do  not converge to second moments of a
Poisson distribution because of "short returns" which will be excluded in our
situation by the assumption (\ref{5.4}) from Section \ref{sec5}. For Harris
recurrent Markov chains convergence in distribution to compound Poisson 
random variables of their number of arrivals
to small sets were studied in \cite{Er} (see
also references there). It seems that many of the results about convergence
to compound Poisson distributions can be extended to the nonconventional 
setup in the spirit of the present paper. This can be done not relying on the
 basic result from 
 \cite{Se} but by proving directly convergence of moments employing 
 combinatorial arguments similar to ones in the proofs below.
 
\end{remark}

\section{A nonconventional Poisson theorem}\label{sec3}
\setcounter{equation}{0}

For any ordered collection of $\ell$ indices $J=(j_1,...,j_\ell)\in\bbZ_+^\ell,
\, j_1<...<j_\ell$ set
\begin{equation}\label{3.1}
X_J=\left\{\begin{array}{ll}
\xi_{q_1(i)}\xi_{q_2(i)}\cdots\xi_{q_\ell(i)}\,\,\mbox{if}\,\, j_l=q_l(i)\,\,
\mbox{for}\,\, i\leq n\,\,\mbox{and all}\,\, l=1,...,\ell\\
0\,\,\mbox{if there is no}\,\, i\geq 1,\, i\leq n\,\,\mbox{such that}\,\,
j_l=q_l(i)\,\,\forall l=1,...,\ell.
\end{array}\right.
\end{equation}
The collection $\{ X_J\}$ is an example of a so called dissociated family
of random variables which means that if $\{ X_J\}_{J\in\cJ}$ and 
$\{ X_K\}_{K\in\cK}$ are two subfamilies and $(\cup_{J\in\cJ}J)\cap
(\cup_{K\in\cK}K)=\emptyset$ then these subfamilies are independent. Thus,
we can apply Theorem 2 from \cite{BE} (see also Section 2.3 and 9.3 in 
\cite{BHJ}) which yields that
\begin{equation}\label{3.2}
\sup_{\Gam\subset\bbZ_+}|P\{ S_n\in \Gam\}-P_{\la_n}(\Gam)|\leq\
min(1,\la^{-1}_n)\big(I_1(n)+I_2(n)+I_3(n)\big)
\end{equation}
where $S_n$ is defined by (\ref{1.5}), $\la_n=np^\ell_n$,
\begin{equation*}
I_1(n)=\sum_Jp^2_J,\, I_2(n)=\sum_J\sum_{K\ne J,\, K\cap J\ne\emptyset}p_Jp_K,
\, I_3(n)=\sum_J\sum_{K\ne J,\, K\cap J\ne\emptyset}EX_JX_K
\end{equation*}
and
\begin{equation}\label{3.3}
p_{(j_1,...,j_\ell)}=\left\{\begin{array}{ll}
p_n^\ell\,\,\mbox{if}\,\, j_l=q_l(i)\,\,\forall l=1,...,\ell\\
0\,\,\mbox{if there is no}\,\, i\geq 1\,\,\mbox{such that}\,\, j_l=q_l(i)\,\,
\forall l=1,...,\ell.
\end{array}\right.
\end{equation}

Clearly,
\begin{equation}\label{3.4}
I_1(n)=np_n^{2\ell}=p_n^\ell\la_n.
\end{equation}
Observe that if $J=(j_1,...,j_\ell),\, p_J\ne 0$ is fixed and 
$K=(k_1,...,k_\ell),\, p_K\ne 0$ satisfies $K\cap J\supset\{ j_{l_1}\}=
\{ k_{l_2}\}$ then such $K$ is uniquely determined by $l_1$ and $l_2$. Hence,
\begin{equation}\label{3.5}
I_2(n)\leq n\ell^2p_n^{2\ell}=\ell^2p_n^\ell\la_n.
\end{equation}
Next, if $K\ne J$ then $EX_JX_K\leq p_n^{\ell+1}$ and since by the above
 argument there are no more than $n\ell^2$ terms in the sum for $I_3(n)$ we
 obtain that
 \begin{equation}\label{3.6}
 I_3(n)\leq n\ell^2p_n^{\ell+1}=\ell^2p_n\la_n.
 \end{equation}
 It follows from (\ref{3.2})--(\ref{3.6}) that
 \begin{equation*}
 |P\{ S_n\in\Gam\}-P_{\la_n}(\Gam)|\leq (2\ell^2+1)p_n
 \end{equation*}
 which yields (\ref{2.2}) taking into account that
 \begin{equation*}
 |P_\la(\Gam)-P_{\la_n}(\Gam)|\leq 2|\la-\la_n|e^{\max(\la,\la_n)}.
 \end{equation*}
 \qed
 
 Next, we formulate the main result from \cite{Se} which we will rely upon
 in Sections \ref{sec4} and \ref{sec5} but also, as a warm up, we will use 
 it below in the
 simpler situation of this section in order to obtain an alternative proof
 of convergence of $P\{ S_n\in A\}$ to $P_\la(A)$ as $n\to\infty$ though
 without an error estimate as in the right hand side of (\ref{2.2}). 
 
 \begin{theorem}\label{thm3.1} (\cite{Se}) Let $\eta_1^{(n)},...,\eta_n^{(n)},
 \, n=1,2,...$ be an array of 0--1 random variables, $J_r(n),\, r\leq n$ be
 the family of all $r$-tuples $(i_1,i_2,...,i_r)$ of mutually distinct indices
  between 1 and $n$ and for any $(i_1,...,i_r)\in J_r(n)$ set
  \[
  b^{(n)}_{i_1,...,i_r}=P\{\eta^{(n)}_{i_1}=...=\eta^{(n)}_{i_r}=1\}.
  \]
  Assume that
  \begin{equation}\label{3.7}
  \lim_{n\to\infty}\max_{1\leq i\leq n}b_i^{(n)}=0,\,\,\lim_{n\to\infty}
  \sum_{i=1}^nb_i^{(n)}=\la>0,
  \end{equation}
  for $n=1,2,...$ there exist "rare" sets $I_r(n)\subset J_r(n)$ such that
  \begin{equation}\label{3.8}
  \lim_{n\to\infty}\sum_{(i_1,...,i_r)\in I_r(n)}b^{(n)}_{i_1...i_r}=
  \lim_{n\to\infty}\sum_{(i_1,...,i_r)\in I_r(n)}b^{(n)}_{i_1}\cdots b_{i_r}=0
  \end{equation}
  and uniformly in $(i_1,...,i_r)\in J_r(n)\setminus I_r(n)$,
  \begin{equation}\label{3.9}
  \lim_{n\to\infty}\frac {b^{(n)}_{i_1...i_r}}{b^{(n)}_{i_1}\cdots 
  b^{(n)}_{i_r}}=1.
  \end{equation}
  Then for $S_n=\sum_{i=1}^n\eta_i^{(n)}$,
  \begin{equation}\label{3.10}
  \lim_{n\to\infty}P\{ S_n=k\}=\frac {\la^ke^{-\la}}{k!},\, k=0,1,2,....
  \end{equation}
  \end{theorem}
  
  Now we set
  \[
  \eta_i^{(n)}=\xi^{(n)}_{q_1(i)}\xi^{(n)}_{q_2(i)}\cdots\xi^{(n)}_{q_1(i)}
  \]
  and check the conditions of Theorem \ref{thm3.1}. For any two positive
  integers $l,\tilde l$ set
  \begin{equation}\label{3.11}
  \rho(l,\tilde l)=\min_{1\leq i,j\leq\ell}|q_i(l)-q_j(\tilde l)|.
  \end{equation}
  A sequence $J=\{ j_1,j_2,...,j_l\}$ of distinct positive integers will be
  called a cluster here if for any $j,\tilde j\in J$ there exists a chain 
  $j_{i_1}=j,j_{i_2},...,j_{i_{m-1}},j_{i_m}=\tilde j$ of integers from $J$
  such that $\rho(j_{i_k},
  j_{i_{k+1}})=0$ for all $k=1,2,...,m-1$. Suppose that $J$ is a part of
  another finite sequence $\tilde J$ of distinct positive integers then
  we say that $J$ is a maximal cluster in $\tilde J$ if $J\cup\{\tilde j\}$
  is already not a cluster for any $\tilde j\in\tilde J$.
  In the notations of Theorem \ref{3.1} we define now "rare" sets $I_r(n)$
  by $I_r(n)=\cup_{1\leq k\leq r-1}I_r^{(k)}(n)$ where $I_r^{(k)}(n)$ is
  the collection of all $r$-tuples from $J_r(n)$ which contain exactly $k$
  maximal clusters . Hence, $J_r(n)\setminus I_r(n)$ consists of $r$-tuples 
  whose all maxinal clusters are singeltons.
  
  Clearly, 
  \[
  b_i^{(n)}=p_i^{(n)}\to 0\,\,\mbox{and}\,\, \sum_{i=1}^nb_i^{(n)}=np_n^\ell
  =\la_n\to\la\,\,\mbox{as}\,\, n\to\infty,
  \]
  and so (\ref{3.7}) holds true. If $(i_1,...,i_r)\in J_r(n)\setminus I_r(n)$
  then $b_{i_1...i_r}=p_n^{\ell r}=b_{i_1}\cdots b_{i_r}$, and so (\ref{3.9})
  is satisfied, as well. Next, if $l$ is fixed and we know that $q_i(l)=q_j(m)$
  then $m$ is uniquely determined by $l,i$ and $j$. Hence, if $l$ is fixed and
  $\rho(l,m)=0$
  then there exist no more than $\ell^2$ possibilities for $m$. It follows 
  that there are no more than $r!n\ell^{2r}$ possibilities for the choice of 
  numbers in any cluster in each sequence from $J_r(n)$. Hence,
  \begin{equation}\label{3.12}
  \#(I_r^{(k)}(n))\leq(r!n\ell^{2r})^k,\, k=1,...,r-1
  \end{equation}
  where $\# \Gam$ denotes the cardinality of a set $\Gam$.
  
   Now, observe that for any sequence $(i_1,...,i_l),\, l\geq 2$ of distinct 
   indices, in particular, for a cluster $b_{i_1...i_l}\leq p_n^{\ell+1}$, and
   so 
   \begin{equation}\label{3.13}
   b_{i_1...i_r}\leq p_n^{k\ell+1}\,\,\mbox{for any}\,\,(i_1,...,i_r)\in 
   I^{(k)}_r(n),\, k\leq r-1
   \end{equation}
   since each $I_r^{(k)}(n)$ with $k\leq r-1$ contains at least one cluster
   which is not a singelton. Hence, by (\ref{3.12}) and (\ref{3.13}),
   \begin{equation}\label{3.14}
   \sum_{(i_1,...,i_r)\in I_r(n)}b_{i_1...i_r}\leq p_n\sum_{k=1}^{r-1}\la_n^k
   (r!\ell^{2r})^k\to 0\,\,\mbox{as}\,\, n\to\infty
   \end{equation}
   while
   \begin{equation}\label{3.15}
   \sum_{(i_1,...,i_r)\in I_r(n)}b_{i_1}\cdots b_{i_r}\leq\sum_{k=1}^{r-1}
   \la_n^kp_n^{(r-k)\ell}(r!\ell^{2r})^k\to 0\,\,\mbox{as}\,\, n\to\infty
   \end{equation}
   implying (\ref{3.8}) and completing the proof of (\ref{3.10}). \qed

 \section{Poisson limits for arrivals to small sets: Markov chains}\label{sec4}
\setcounter{equation}{0}

Next, we prove Theorem \ref{thm2.2}. Set
\[
a(l)=\min(\ln l,\,\min_{1\leq i\leq\ell -1}(q_{i+1}(l)-q_i(l)).
\]
Now a sequence $J=\{ j_1,j_2,...,j_l\}$ of distinct positive integers will
be called an $(a,n)$-cluster if for any $j,\tilde j\in J$ there exists a
chain $j_{i_1}=j,j_{i_2},...,j_{i_{m-1}},j_{i_m}=\tilde j$ of integers
from $J$ such that 
\begin{equation}\label{4.1}
\rho(j_{i_k},j_{i_{k+1}})\leq a(n)\,\,\,\forall k=1,2,...,m-1
\end{equation}
with $\rho(\cdot,\cdot)$ defined by (\ref{3.11}). The definition of
maximal clusters remains as before while we define rare sets $I_r(n)$ as
collections of $r$-tuples $J=(i_1,i_2,...,i_r)$ from $J_r(n)$ 
(which, recall, denotes the collection of all $r$-tuples of mutually distinct
indices between 1 and $n$) which either contain a
cluster containing more than one element or
\begin{equation}\label{4.2}
i_{\mbox{min}}(J)=\min_{1\leq l\leq r}i_l\leq a(n).
\end{equation}

By (\ref{2.3}) and the Radon--Nikodim theorem there exists a transition 
density $p(x,y)$ so that 
\begin{equation*}
P(x,\Gam)=\int_\Gam p(x,y)dm(y)\,\,\mbox{and}\,\, C^{-1}\leq p(x,y)\leq C.
\end{equation*}
Let $p(n,x,y)$ denotes the $n$-step transition density so that
\[
P(n,x,\Gam)=P_x\{ X_n\in\Gam\}=P\{ X_n\in\Gam|X_0=x\}=\int_\Gam p(n,x,y)dm(y).
\]
We will rely on the well known fact (see, for instance, \cite{Doo}, \S 5 
in Ch. V) that under (\ref{2.3}) there exists a unique invariant measure 
$\mu$ (i.e. $\int_Md\mu(x)P(x,\Gam)=\mu(\Gam),\,\forall\Gam\in\cB$) having
a density $p(x)=\frac {d\mu(x)}{dm(x)}$ satisfying $C^{-1}\leq p(x)\leq C$
and
\begin{equation}\label{4.3}
\sup_{x,y\in M}|p(n,x,y)-p(y)|\leq C_1e^{-\be n}\,\,\forall n\geq 1
\end{equation}
for some $C_1,\be>0$ independent of $n$.

Set
\[
\eta_l^{(n)}=\prod_{j=1}^\ell\bbI_{\Gam_n}(X_{q_j(l)}),\, l=1,2,....
\]
We start verifying conditions of Theorem \ref{thm3.1} in our situation 
observing that by the Chapman--Kolmogorov formula
\begin{eqnarray}\label{4.4}
&\quad b_l^{(n)}=P_\nu\{\eta_l^{(n)}=1\}=\int_Md\nu(x)\int_{\Gam_n}p(q_1(l),
x,x_1)\int_{\Gam_n}p(q_2(l)-q_1(l),x_1,x_2)\\
&...\int_{\Gam_n}p(q_\ell(l)-q_{\ell-1}(l),x_{\ell-1},x_\ell)
dm(x_1)...dm(x_\ell)\nonumber
\end{eqnarray}
where, recall, $\nu$ is the initial distribution of the Markov chain. Since
$q_1(l)\geq l$ we obtain from (\ref{2.4}), (\ref{4.3}) and (\ref{4.4})
that for any $l\geq 1$,
\begin{equation}\label{4.5}
|b_l^{(n)}-(\mu(\Gam_n))^\ell|\leq C_2(\mu(\Gam_n))^\ell\exp(-\be a(l))
\end{equation}
for some $C_2>0$ independent of $l$. Hence, by (\ref{1.7}) and (\ref{4.5}),
\begin{equation}\label{4.6}
\max_{1\leq l\leq n}b_l^{(n)}\leq (C_2+1)(\mu(\Gam_n))^\ell\to 0\,\,\mbox{as}
\,\,n\to\infty
\end{equation}
and
\begin{equation}\label{4.7}
|\sum_{l=1}^nb_l^{(n)}-n(\mu(\Gam_n))^\ell|\leq C_2((\mu(\Gam_n))^\ell
\sum_{l=1}^ne^{-\be a(l)}\to 0\,\,\mbox{as}\,\, n\to\infty,
\end{equation}
and so (\ref{3.7}) is satisfied.

Observe that since
$p(x,y)\leq C$ then $p(l,x,y)\leq C$ for all $l,x,y$ in view of the
Chapman--Kolmogorov formula. Let $(i_1,...,i_l)$ be a sequence of
distinct integers such that for some pairs $(m_1,i_{j_1}),...,(m_k,i_{j_k})$,
\begin{equation}\label{4.8}
q_{m_1}(i_{j_1})<q_{m_2}(i_{j_2})<...<q_{m_k}(i_{j_k})
\end{equation}
where pairs are different but either $i$'s or $m$'s may repeat themselves.
It follows by the Chapman--Kolmogorov formula that
\begin{eqnarray}\label{4.9}
&\quad b^{(n)}_{i_1...i_l}=P_\nu\{\eta_{i_1}^{(n)}=1,...,\eta_{i_l}^{(n)}
=1\}\\
&\leq P_\nu\{ X_{q_{m_1}}(i_{j_1})\in\Gam_n,..., X_{q_{m_k}}(i_{j_k})
\in\Gam_n\}\nonumber\\
&=\int_Md\nu(x)\int_{\Gam_n}p(q_{m_1}(i_{j_1}),x,x_1)\int_{\Gam_n}
p(q_{m_2}(i_{j_2})-q_{m_1}(i_{j_1}),x_1,x_2)\nonumber\\
&...\int_{\Gam_n}p(q_{m_k}(i_{j_k})-q_{m_{k-1}}(i_{j_{k-1}}),x_{k-1},x_k)
dm(x_1)...dm(x_k)\nonumber\\
&\leq(Cm(\Gam_n))^k\leq(C^2\mu(\Gam_n))^k.\nonumber
\end{eqnarray}

Next, we represent the rare sets $I_r(n)$ in the form
\begin{equation}\label{4.10}
I_r(n)=\big(\cup_{1\leq k\leq r}I^{(k,1)}_r(n)\big)\cup
\big(\cup_{k=1}^{r-1}I^{k,0}_r(n)\big)
\end{equation}
where $I^{(k,1)}_r(n)$ and $I^{(k,1)}_r(n)$ are the sets of $r$-tuples from
$J_r(n),\, r\leq n$ which contain exactly $k$ maximal $(a,n)$-clusters and
each $I^{(k,1)}_r(n)$ contains an $(a,n)$-cluster $J$ with 
$i_{\mbox{min}}(J)\leq a(n)$ while $I^{(k,0)}_r(n)$ contains no such
$(a,n)$-clusters. Clearly, no $r$-tuple from $I_r(n)$ may contain more than
 one maximal $(a,n)$-cluster $J$ satisfying $i_{\mbox{min}}(J)\leq a(n)$.
Observe that
 $I^{(r,l)}_r(n),\, l=0,1$ contains only singelton $(a,n)$-clusters, and so 
it is a subset of $I_r(n)$ only if $l=1$.  In order to
estimate the cardinality of rare sets observe that if $l$ and $m$ belong
to a same $(a,n)$-cluster and $l$ is fixed then there are no more than
$2^r(a(n))^r\ell^{r+1}$ possibilities for $m$. It follows that there are no more
than $n(a(n))^{r^2}2^{r^2}\ell^{r(r+1)}r!$ possibilities for the choice of numbers in any
cluster in $J=(i_1,...,i_l)$ with $i_{\mbox{min}}(J)>a(n)$ while if
$i_{\mbox{min}}(J)\leq a(n)$ then there are no more than 
$(a(n))^{r^2+1}2^{r^2}\ell^{r(r+1)}r!$ such choices. Hence,
\begin{equation}\label{4.11}
\#(I_r^{(k,l)}(n))\leq (2^{r^2}\ell^{r(r+1)}r!)^k(a(n))^{kr^2+l}n^{k-l},\,\, 
k=1,2,...,r; l=0,1.
\end{equation}
where $\# \Gam$ denotes cardinality of a set $\Gam$. 

Now observe that each cluster $J=(i_1,...,i_l)$ which is not a singleton
yields at least $\ell+1$ pairs $(m_1,i_{j_1}),...,(m_{\ell+1},i_{j_{\ell+1}})$
satisfying (\ref{4.8}). Hence, any $\tilde J=(\tilde i_1,
...,\tilde i_r)\in I_r^{(k,l)}(n)$ with $k<r$ yields at least $k\ell+1$
pairs satisfying (\ref{4.8})
while each $r$-tuple from $I_r^{(r,l)}(n)$ yields at least $r\ell$ such
pairs. Since $a(n)\leq\ln n$ then these arguments together with (\ref{1.7})
and (\ref{4.9})--(\ref{4.11}) yield that
\begin{eqnarray}\label{4.12}
&\sum_{(i_1,...,i_r)\in I_r(n)}b^{(n)}_{i_1...i_r}\leq \sum_{k=1}^r
\sum_{(i_1,...,i_r)\in I_r^{(k,1)}(n)}b^{(n)}_{i_1...i_r}\\
&+\sum_{k=1}^{r-1}\sum_{(i_1,...,i_r)\in I_r^{(k,0)}(n)}b^{(n)}_{i_1...i_r}
\leq C_3\big( \sum_{k=1}^r((a(n))^{r^2k+1}(\mu(\Gam_n))^\ell\la_n^{k-1}
\nonumber\\
&+\sum_{k=1}^{r-1}(a(n))^{r^2k+1}\la_n^k\mu(\Gam_n)\big)\to 0\,\,\mbox{as}
\,\, n\to\infty\nonumber
\end{eqnarray}
and
\begin{eqnarray}\label{4.13}
&\sum_{(i_1,...,i_r)\in I_r(n)}b^{(n)}_{i_1}\cdots b^{(n)}_{i_r}\\
&\leq \sum_{k=1}^r\sum_{(i_1,...,i_r)\in I_r^{(k,1)}(n)}
b^{(n)}_{i_1}\cdots b^{(n)}_{i_r}\nonumber\\
&+\sum_{k=1}^{r-1}\sum_{(i_1,...,i_r)\in I_r^{(k,0)}(n)}b^{(n)}_{i_1}\cdots
b^{(n)}_{i_r}\nonumber\\
&\leq C_3\big( \sum_{k=1}^r((a(n))^{r^2k+1}(\mu(\Gam_n))^\ell\la_n^{k-1}
\nonumber\\
&+\sum_{k=1}^{r-1}(a(n))^{r^2k+1}\la_n^k(\mu(\Gam_n))^{(r-k)\ell}\to 0\,\,\mbox{as}
\,\, n\to\infty\nonumber
\end{eqnarray}
where $\la_n=n(\mu(\Gam_n))^\ell$ and $C_3>0$ does not depend on $n$, which 
gives (\ref{3.8}).

Next, let $(i_1,...,i_r)\in J_r(n)\setminus I_r(n)$. Then there exist pairs
$(m_1,i_{j_1}),\,(m_2,i_{j_2}),...,(m_{r\ell},i_{j_{r\ell}})$ such that
\begin{equation}\label{4.14}
i_{j_1}\geq a(n)\,\,\,\mbox{and}\,\,\, |q_{m_{l+1}}(i_{j_{l+1}})-q_{m_l}
(i_{j_l})|\geq a(n)\,\,\mbox{for}\,\, l=1,2,...,r\ell-1.
\end{equation}
Employing again the Chapman--Kolmogorov formula together with (\ref{4.3})
similarly to (\ref{4.5}) we obtain that for such $(i_1,...,i_r)$,
\begin{equation}\label{4.15}
|b^{(n)}_{i_1...i_r}-(\mu(\Gam_n))^{r\ell}|\leq C_4\exp(-\be a(n))
(\mu(\Gam_n))^{r\ell}
\end{equation}
for some $C_4>0$ independent of $n$. This together with (\ref{4.5}) yields
(\ref{3.9}) and completes the proof of Theorem \ref{thm2.2}.  \qed

\section{Nonconventional Poisson limits for subshifts}\label{sec5}
\setcounter{equation}{0}

Before proving Theorem \ref{thm2.3} itself we recall few basic facts 
about Gibbs measures for topologically mixing subshifts of finite type
whose proofs can be found, for instance, in Ch.1 of \cite{Bo}. Namely,
in addition to (\ref{2.6}) we have that for $P$-almost all $\om$,
\begin{equation}\label{5.1}
\lim_{n\to\infty}\frac 1n\sum_{i=0}^n\phi\circ T^i(\om)=\int_\Om\phi dP=
\Pi-h_P(T)
\end{equation}
where $h_P(T)>0$ is the Kolmogorov--Sinai entropy of $T$ with respect to
$P$. It follows from (\ref{2.6}) and (\ref{5.1}) that for $P$-almost all
$\om=(\om_0,\om_1,...)$,
\begin{equation}\label{5.2}
\lim_{n\to\infty}\frac 1n\ln P([\om_0,\om_1...,\om_{n-1}])=-h_P(T)<0.
\end{equation}
Another important fact which we need is the exponentially fast $\psi$-mixing
of such subshifts, namely, that there exist constants $\be,C>0$ such that
for any two cylinder sets $U=[a_0,a_1,...,a_l]\subset\Om$ and $V=[b_0,b_1,...,
b_m]\subset\Om$,
\begin{equation}\label{5.3}
|P(U\cap T^{-n}V)-P(U)P(V)|\leq Ce^{-\be(n-l)}P(U)P(V)
\end{equation}
provided $n\geq l+1$.

Set $a(n)=[\ln^{1+\ve}n],\,\ve\in(0,\gam)$ where $\gam$ is the same as in
(\ref{2.7}). Denote by $\Om^*$ the set of all $\om^*\in\Om$ such that 
\begin{equation}\label{5.4}
C_n(\om^*)\cap T^{-i}C_n(\om^*)=\emptyset\quad\mbox{for all}\quad i=1,2,...,
a(n)
\end{equation}
and
\begin{equation}\label{5.5}
\lim_{n\to\infty}\frac 1n\sum_{i=0}^{n-1}\phi\circ T^i(\om_n)=\Pi-h_P(T)
\quad\mbox{for any sequence}\quad \om_n\in C_n(\om^*).
\end{equation}
 It is clear that 
$\Om^*$ is measurable and it was shown in \cite{Pi} and \cite{De} that
$P(\Om^*)=1$. For reader's convenience we recall the corresponding argument
observing, first, that (\ref{5.5}) follows from (\ref{5.1}). As to (\ref{5.4})
we note that $C_n\cap T^{-i}C_n\ne\emptyset$ for a cylinder set $C_n$ if and 
only if it contains a periodic point of period $i$. The number of periodic
points of period $l$ grows with $l$ exponentially, and so there are at most 
$g^{a(n)}$ periodic points of periods up to $a(n)$ for some $g>1$. Thus,
there exist at most $g^{a(n)}$ cylinder sets $C_n$ of length $n$ such that
$C_n\cap T^{-i}C_n\ne\emptyset$ for some $i\leq a(n)$. On the other hand,
relying on (\ref{5.2}) we see that $P(C_n(\om^*))$, and so also 
$g^{a(n)}P(C_n(\om^*))$, decay exponentially fast in $n$ for $P$-almost
all $\om^*$. This together with the Borel--Cantelli lemma completes the
argument.

Next, we return to our nonconventional setup. Put
\[
\eta_l^{(n)}(\om)=\prod_{j=1}^\ell\bbI_{B_n(\om^*)}\circ T^{q_j(l)}(\om),
\,\, l=1,2,...,N_n(\om^*)
\]
where $\om^*\in\Om^*$, $N_n(\om^*)$ was defined in Theorem \ref{thm2.3},
$B_n(\om^*)\subset C_n(\om^*)$ and $B_n(\om^*)\in\cF_{n+[s\ln n]},\, 
n=1,2,...$. Observe a slight change of notations here in comparison to
Theorem \ref{thm3.1} by writing $\eta^{(n)}_l$ in place of 
$\eta^{(N_n(\om^*))}$ which would be unwieldy. We start verifying the 
conditions of Theorem \ref{thm3.1} observing that
\begin{equation}\label{5.6}
b_l^{(n)}=P\{\eta_l^{(n)}=1\}=P\big(\cap_{j=1}^\ell T^{(q_j(l)-q_1(l))}
B_n(\om^*)\big).
\end{equation}
Taking into account (\ref{2.7}) and applying (\ref{5.3}) repeatedly we
obtain that 
\begin{equation}\label{5.7}
\big\vert b_l^{(n)}-\big(P(B_n(\om^*))\big)^\ell\big\vert\leq D_1
\exp(-\be a(n))\big(P(B_n(\om^*))\big)^\ell
\end{equation}
for all $l\geq L(n)$ where
\begin{equation}\label{5.8}
L(n)=\min\{ k:\, c(\ln k)^{1+\gam}>2(n+a(n))\}=\big[\exp\big(
\frac {2(n+a(n))}{c}\big)^{\frac 1{1+\gam}}\big]+1
\end{equation}
grows in $n$ subexponentially and $D_1>0$ does not depend on $l,n$ and
 $\om^*$.
 
 Since $P(B_n(\om^*))$ decays in $n$ exponentially fast we obtain from
 (\ref{2.8}) and (\ref{5.6})--(\ref{5.8}) that 
 \begin{equation}\label{5.10}
 \max_{1\leq l\leq N_n(\om^*)}b_l^{(n)}\leq P(B_n(\om^*))\to 0\,\,\mbox{as}
 \,\, n\to\infty
 \end{equation}
 and
 \begin{eqnarray}\label{5.11}
 &\big\vert\sum_{l=1}^{N_n(\om^*)}b_l^{(n)}-N_n(\om^*)
 \big(P(B_n(\om^*))\big)^\ell\big\vert\\
 &\leq 2L(N)P(B_n(\om^*))+D_1\big(P(B_n(\om^*))\big)^\ell N_n(\om^*)
 e^{-\be a(n)}\to 0\,\,\mbox{as}\,\, n\to\infty\nonumber
 \end{eqnarray}
 yielding (\ref{3.7}).

We introduce again an $(a,n)$-cluster which is a sequence $J=\{ j_1,j_2,
...,j_l\}$ of distinct positive integers such that for any $j,\tilde j\in J$ 
there exists a chain $j_{i_1}=j,j_{i_2},...,j_{i_{m-1}},j_{i_m}=\tilde j$
of integers from $J$ such that
\begin{equation}\label{5.12}
\rho(j_{i_k},j_{i_{k+1}})\leq n+a(n)\quad\forall k=1,2,...,m-1
\end{equation}
with $\rho(\cdot,\cdot)$ defined by (\ref{3.11}). Denote by $J_r(n)=
J_r(n,\om^*)$ the set of all $r$-tuples $(i_1,...,i_r)$ of distinct integers
between 1 and $N_n(\om^*)$. The definition of maximal clusters remains the same
as before
 and we define rare sets $I_r(n)=I_r(n,\om^*)$ as collections of $r$-tuples
 $J=(i_1,i_2,...,i_r)\in J_r(n)$ which either contain not only singelton
 clusters or
 \begin{equation}\label{5.13}
 i_{\mbox{min}}(J)=\min_{1\leq j\leq r}i_j\leq L(n).
 \end{equation}
 
 Next, we represent the rare sets $I_r(n)$ in the form
\begin{equation}\label{5.14}
I_r(n)=\big(\cup_{1\leq l\leq k\leq r}I^{(k,l)}_r(n)\big)\cup
\big(\cup_{k=1}^{r-1}I^{k,0}_r(n)\big)
\end{equation}
where each $I^{(k,l)}_r(n)\subset J_r(n)$ contains exactly $k$ maximal 
$(a,n)$-clusters while $l$ of them are collections $J$ satisfying 
$i_{\mbox{min}}(J)\leq L(n)$. In order to estimate cardinality of 
$I^{(k,l)}_r(n)$ observe that if $i_1$ and $i_2$ belong to the same
$(a,n)$-cluster and $i_1$ is fixed then there are no more than 
$2^r\ell^{r+1}(n+a(n))^r$ possibilities for $i_2$. If $i\leq L(n)$ then, of course,
$L(n)$ bounds the number of choices for $i$. It follows that there are no more
than $2^{r^2}N_n(\om^*)\ell^{r(r+1)}r!(n+a(n))^{r^2}$ possibilities for the 
choice of 
integers in any cluster in $J=(i_1,...,i_l)$ with $i_{\mbox{min}}(J)>L(n)$ 
while if $i_{\mbox{min}}(J)\leq L(n)$ then there are no more than 
$2^{r^2}L(n)\ell^{r(r+1)}r!(n+a(n))^{r^2}$ such choices. Thus, (\ref{4.11})
 takes here the form
\begin{equation}\label{5.15}
\#(I_r^{(k,l)}(n))\leq (2^{r^2}\ell^{r(r+1)}r!)^k(n+a(n))^{kr^2}(L(n))^l
(N_n(\om^*))^{k-l}
\end{equation}
for $k=1,2,...,r$ and $l=0,1,...,k$.

Let $(i_1,...,i_l)$ be a sequence of
distinct integers such that for some pairs $(m_1,i_{j_1}),...,(m_k,i_{j_k})$,
\begin{eqnarray}\label{5.16}
&q_{m_1}(i_{j_1})<q_{m_2}(i_{j_2})-(n+a(n))<q_{m_3}(i_{j_3})-2(n+a(n))\\
&<...<q_{m_k}(i_{j_k})-(k-1)(n+a(n))\nonumber
\end{eqnarray}
where pairs are different but either $i$'s or $m$'s may repeat themselves.
Then applying repeatedly (\ref{5.3}) we obtain that for all $n$ large enough
(say, when $a(n)>s\ln n$),
\begin{eqnarray}\label{5.17}
&b^{(n)}_{i_1...i_l}\leq P\big( B_n(\om^*)\cap T^{-(q_{m_2}(i_{j_2})-
q_{m_1}(i_{j_1}))}B_n(\om^*)\cap\\
&...\cap T^{-(q_{m_k}(i_{j_k})-q_{m_1}(i_{j_1}))}B_n(\om^*)\big)\leq
D_2\big( P\big( B_n(\om^*)\big)^k\nonumber
\end{eqnarray}
for some $D_2>0$ independent of $n,\om^*,i_1,...,i_l$ and $B_n(\om^*)$.

By the definition each $J\in I_r^{(k,l)}(n)$ consists of
$k$ maximal $(a,n)$-clusters $J_1,...,J_l,...,J_k$ such that 
$i_{\mbox{min}}(J_j)\leq L(n)$ for $j=1,...,l$ and $i_{\mbox{min}}(J_j)
>L(n)$ for $j=l+1,...,k$. For each $j=1,...,k$ choose $\tilde i_j\in J_j$
arbitrarily and set $i_j=\tilde i_j$ for $j=1,...,l$ and 
$i_{l+(j-l)\ell+b}=\tilde i_j$ for $j=l,l+1,...,k-1$ and $b=1,...,\ell$. 
Set also $m_j=1$ for $j=1,...,l$ and $m_{l+(j-l)\ell+b}=b$ for $j=l,l+1,...,
k-1$ and $b=1,...,\ell$. Next, we reorder the pairs $(m_j,i_j)$ so that 
$q_{m_{j_b}}(i_{j_b})$ increases in $b$ and it follows from the definition
of $L(N)$ and of $(a,n)$-clusters that $l+(k-l)\ell$ pairs $(m_{j_b},i_{j_b}),
\, b=1,...,l+(k-l)\ell$ will satisfy (\ref{5.16}). This together with 
(\ref{5.17}) yields that for any $(i_1,...,i_r)\in I_r^{(k,l)}(n)$ and
$n$ large enough,
 \begin{equation}\label{5.18}
 b^{(n)}_{i_1...i_l}\leq D_3\big(P(B_n(\om^*))\big)^{(k-l)\ell+l}
 \end{equation}
 for some $D_3>0$ independent on $n,\om^*$ and $B_n(\om^*)$.
 
 The estimate (\ref{5.18}) will suffice for our purposes when $l\geq 1$
 but for $J=(\tilde i_1,...,\tilde i_r)\in I_r^{(k,0)}(n),\, k<r$ a better 
 estimate of $b^{(n)}_{\tilde i_1...\tilde i_l}$ will be needed. So let
 such $J$ consists of $k$ maximal $(a,n)$-clusters $J_1,...,J_k$ and since
 $k<r$ one of them must be not a singelton. Suppose, for instance, $i,\tilde
 i\in J_1,\, i\ne\tilde i$ and without loss of generality assume that
 \[
 |q_j(i)-q_{\tilde j}(\tilde i)|\leq n+a(n)\,\,\mbox{for some}\,\, j,\tilde j=
 1,...,\ell.
 \]
 Set $i_1=i$ and choose $i_b\in J_b,\, b=2,3,...,k$ arbitrarily. Now, order
 $k\ell +1$ pairs $(\tilde j,\tilde i)$ and $(l, i_b),\, l=1,...,\ell;\,
 b=1,...,k$ to obtain pairs $(m_1,i_{j_1}),\, (m_2,i_{j_2}),...,(m_{k\ell+1},
 i_{j_{k\ell+1}})$ so that $q_{m_b}(i_{j_b})$ is nondecreasing in $b$. Let
 $q_{m_l}(i_{j_l})=q_j(i)$ and assume without loss of generality that
 $q_{\tilde j}(\tilde i)\geq q_j(i)$. Then we must have 
 $q_{\tilde j}(\tilde i)=q_{m_{l+1}}(i_{j_{l+1}})$ and
 \begin{eqnarray*}
 &q_{m_1}(i_{j_1})<q_{m_2}(i_{j_2})-(n+a(n))<...<q_{m_l}(i_{j_l})-
 (l-1)(n+a(n))\\
 &\leq q_{m_{l+1}}(i_{j_{l+1}})-(l-1)(n+a(n))\\
 &<q_{m_{l+2}}(i_{j_{l+2}})-l(n+a(n))<...<q_{m_{k\ell+1}}(i_{j_{k\ell+1}})
 -(k\ell-1)(n+a(n)).
 \end{eqnarray*}
 Applying repeatedly (\ref{5.3}) we obtain from here similarly to (\ref{5.17})
 that for all $n$ large enough, 
\begin{eqnarray}\label{5.19}
& b_{\tilde i_1...\tilde i_r}\leq P\big( B_n(\om^*)\cap T^{q_{m_2}(i_{j_2})
-q_{m_1}(i_{j_1})}B_n(\om^*)\cap ...\\
&...\cap T^{q_{m_{k\ell}}(i_{j_{k\ell}})-q_{m_1}(i_{j_1})}B_n(\om^*)\big)
\nonumber\\
&\leq D_4\big(P(B_n(\om^*))\big)^{k\ell-1}P\big( B_n(\om^*)\cap 
T^{q_{m_{l+1}}(i_{j_{l+1}})-q_{m_l}(i_{j_l})}B_n(\om^*)\big)\nonumber
\end{eqnarray} 
 for some $D_4>0$ independent of $n,\,\om^*$ and $B_n(\om^*)$. 
 
 By (\ref{5.4}) if
 \[
 q_{m_{l+1}}(i_{j_{l+1}})-q_{m_l}(i_{j_l})\leq a(n)
 \]
 then
 \[
 B_n(\om^*)\cap T^{q_{m_{l+1}}(i_{j_{l+1}})-q_{m_l}(i_{j_l})}B_n(\om^*)=
 \emptyset,
 \]
 and so the right hand side of (\ref{5.19}) is zero. If
 \[
 n+[s\ln n]<q_{m_{l+1}}(i_{j_{l+1}})-q_{m_l}(i_{j_l})\leq n+a(n)
 \]
 then we still can use (\ref{5.3}) to obtain that
 \begin{equation}\label{5.20}
 P\big( B_n(\om^*)\cap T^{q_{m_{l+1}}(i_{j_{l+1}})-q_{m_l}(i_{j_l})}
 B_n(\om^*)\big)\leq C\big(P(B_n(\om^*))\big)^2.
 \end{equation}
 Now let
 \begin{equation}\label{5.21}
a(n) <q_{m_{l+1}}(i_{j_{l+1}})-q_{m_l}(i_{j_l})\leq n+[s\ln n].
 \end{equation}
 We can represent $B_n(\om^*)$ as a disjoint union
 \[
 B_n(\om^*)=\cup_{b=1}^{m_n(\om^*)}C^{(b)}_{n,s}(\om^*)
 \]
 where $C^{(b)}_{n,s}(\om^*),\, b=1,...,m_n(\om^*)$ are cylinder sets of
 the length $n+[s\ln n]$ and $m_n(\om^*)\leq n^{\eta s}$ for some $\eta>0$. 
 Employing (\ref{2.6}) and (\ref{5.5}) we obtain that under (\ref{5.21}) for
 any $b,d=1,...,m_n(\om^*)$,
  \begin{equation}\label{5.22}
  P\big( C^{(b)}_{n,s}(\om^*)\cap T^{q_{m_{l+1}}(i_{j_{l+1}})-
  q_{m_l}(i_{j_l})}C^{(d)}_{n,s}(\om^*)\big)\leq D_5P(C^{(b)}_{n,s}(\om^*))
  \exp(-\del a(n))
 \end{equation}
 for some $D_5,\del>0$ where we can take $\del=\frac 12h_P(T)$ provided 
 $n$ is large enough. Summing in $b$ and $d$ in (\ref{5.22}) we obtain that
 \begin{equation}\label{5.23}
 P\big( B_n(\om^*)\cap T^{q_{m_{l+1}}(i_{j_{l+1}})-q_{m_l}(i_{j_l})}
 B_n(\om^*)\big)\leq D_5n^{\eta s}P(B_n(\om^*))\exp(-\del a(n))
 \end{equation}
 provided (\ref{5.21}) holds true. In view of (\ref{2.6}) and (\ref{5.2})
 the estimate (\ref{5.23}) is weaker than (\ref{5.20}) so we will use the former
 in both cases. Hence, from (\ref{5.19}), (\ref{5.20}) and (\ref{5.23}) it
 follows that for any $(\tilde i_1,...,\tilde i_r)\in I^{(k,0)}_r(n)$ for
 all $n$ large enough
 \begin{equation}\label{5.24}
 b^{(n)}_{\tilde i_1,...,\tilde i_r}\leq D_4D_5n^{\eta s}\big(P(B_n(\om^*))
 \big)^{k\ell}\exp(-\del a(n)).
 \end{equation}
  
Now (\ref{5.10}), (\ref{5.15}), (\ref{5.18}), (\ref{5.24}) together with the
definition of $a(n)$ yield that
\begin{eqnarray}\label{5.25}
&\sum_{(i_1,...,i_r)\in I_r(n)}b^{(n)}_{i_1...i_r}\leq \sum_{k=1}^r\sum_{l=1}^k
\sum_{(i_1,...,i_r)\in I_r^{(k,l)}(n)}b^{(n)}_{i_1...i_r}\\
&+\sum_{k=1}^{r-1}\sum_{(i_1,...,i_r)\in I_r^{(k,0)}(n)}b^{(n)}_{i_1...i_r}
\nonumber\\
&\leq D_6\big( \sum_{k=1}^r(n+a(n))^{r^2k}(L(n))^l\big(P(B_n(\om^*))\big)^l
(\la_n(\om^*))^{k-l}
\nonumber\\
&+\sum_{k=1}^{r-1}(n+a(n))^{r^2k}(\la_n(\om^*))^kn^{\eta s}\exp(-\del a(n))\big)
\to 0\,\,\mbox{as}\,\, n\to\infty\nonumber
\end{eqnarray}
and
\begin{eqnarray}\label{5.26}
&\sum_{(i_1,...,i_r)\in I_r(n)}b^{(n)}_{i_1}\cdots b^{(n)}_{i_r}\\
&\leq \sum_{k=1}^r\sum_{l=1}^k\sum_{(i_1,...,i_r)\in I_r^{(k,l)}(n)}
b^{(n)}_{i_1}\cdots b^{(n)}_{i_r}\nonumber\\
&+\sum_{k=1}^{r-1}\sum_{(i_1,...,i_r)\in I_r^{(k,0)}(n)}b^{(n)}_{i_1}\cdots
b^{(n)}_{i_r}\nonumber\\
&\leq D_6\big( \sum_{k=1}^r\sum_{l=1}^k((n+a(n))^{r^2k}(L(n))^\ell
(\la_n(\om^*))^{k-l}(P(B_n(\om^*))^l
\nonumber\\
&+\sum_{k=1}^{r-1}(n+a(n))^{r^2k}(\la_n(\om^*))^k(P(B_n(\om^*)))^{(r-k)\ell}
\to 0\,\,\mbox{as}\,\, n\to\infty\nonumber
\end{eqnarray}
where $\la_n(\om^*)=N_n(\om^*)(P(B_n(\om^*)))^{(r-k)\ell}$ and $D_6>0$ does 
not depend on $n$, which gives (\ref{3.8}).

Next, let $(i_1,...,i_r)\in J_r(n)\setminus I_r(n)$. Then there exist pairs
$(m_1,i_{j_1}),\,(m_2,i_{j_2}),...,(m_{r\ell},i_{j_{r\ell}})$ such that
\begin{equation}\label{5.27}
i_{j_1}\geq L(n)\,\,\,\mbox{and}\,\,\, q_{m_{l+1}}(i_{j_{l+1}})-q_{m_l}
(i_{j_l})> n+a(n)\,\,\mbox{for}\,\, l=1,2,...,r\ell-1.
\end{equation}
Then applying repeatedly (\ref{5.3}) we obtain similarly to (\ref{5.7})
that for such $(i_1,...,i_r)$,
\begin{equation}\label{5.28}
|b^{(n)}_{i_1...i_r}-\big(P(B_n(\om^*))\big)^{r\ell}|\leq D_7\exp(-\be (a(n)-
[s\ln n]))\big(P(B_n(\om^*))\big)^{r\ell}
\end{equation}
for some $D_7>0$ independent of $n,\om^*$ and $B_n(\om^*)$. This together with
(\ref{5.7}) yields (\ref{3.9}) and completes the proof of Theorem \ref{thm2.3}.  \qed

 \bibliography{matz_nonarticles,matz_articles}
\bibliographystyle{alpha}

\end{document}